\documentclass[preprint,prd,aps,showpacs,showkeys,onecolumn,english]{revtex4}
\usepackage{amsmath,amssymb} 
\usepackage[dvips]{graphicx} 
\usepackage{graphics}
\usepackage{graphicx}
\usepackage{epsfig}
\usepackage{slashed}
\usepackage[english]{babel}
\usepackage{amsmath}
\usepackage{amssymb}
\usepackage{graphics}
\usepackage{graphicx}
\usepackage{epsfig}
\usepackage{slashed}
\usepackage{subfigure}	
\usepackage[export]{adjustbox}

\newcommand\nn{\nonumber}
\newcommand\ba{\begin{eqnarray}}
\newcommand\ea{\end{eqnarray}}

\begin{document}

\title{Existence and uniqueness of the solution of a mixed problem for a parabolic equation under nonconventional boundary conditions}
\author{Yu.A. Mammadov\footnote{E-mail: yusifmamedov@icloud.com }, H.I. Ahmadov\footnote{E-mail: hikmatahmadov@yahoo.com}}
\affiliation{Department of Equations of Mathematical Physics,
Faculty of Applied Mathematics and Cybernetics, Baku State
University, st. Z. Khalilov 23, AZ-1148, Baku, Azerbaijan}

\date{\today}

\begin{abstract}
In this study, we investigate a mixed problem linked to a second-order parabolic equation,
characterized by temporal dependencies and variable~coefficients, and constrained
by non-local, non-self-adjoint boundary conditions.
By defining precise conditions on the input data, we establish the unique
solvability of the problem through a synthesis of the residue and
contour integral methods. Moreover, our research yields an
explicit analytical solution, facilitating the direct resolution
of the stated problem.

\end{abstract}

\maketitle

\section{Introduction}
\label{Introduction}
\bigskip

Nonlocal boundary conditions, pioneered by A.A. Samarskiy and A.V.
Bitsadze, play a crucial role in the theory of differential
equations and mathematical physics \cite{1}. These conditions are
integral to the study of multipoint boundary value problems for
ordinary differential equations, which~are extensively applied~in
diverse fields such as electric power grids, electric railway
systems, telecommunication lines, and the chemical analysis of
kinetic reactions \cite{2,3,4,5,6,7,8}.

Despite the broad study of such applications, less attention has
been given~to non-stationary problems incorporating multipoint
boundary conditions. Notably, a specific subset of this research
focuses on three-point boundary conditions within the context of
nonlinear parabolic Cauchy problems, while others explore the
semi-linear parabolic equations in one dimension, introducing
conditions that link boundary points to interior points of the
domain, thus forming a four-point boundary condition framework
\cite{9,10,11,12,13,14,15,16,17}.

Innovative regularization methods for boundary value problems in
parabolic equations have been developed, addressing issues such as
singular perturbations on semi-axes with rational turning
points,~and applying principles like the maximum principle to
verify the asymptotic convergence of solutions \cite{15}.
The examination of~initial boundary value problems for quasi-linear
pseudo-parabolic equations of fractional order also features
prominently, where the Galerkin method facilitates the
establishment of weak solution existences. This research leverages
Sobolev embedding theorems to derive a priori estimates, employing
these along with the Rellich--Kondrashov theorem to substantiate
the existence of solutions under conditions defined by nonlinear,
fractional differentiation operators in the boundary conditions
\cite{16}. An inverse problem for time fractional parabolic partial
differential equations with nonlocal boundary conditions has
also~been addressed. This research utilizes Dirichlet-measured
output data to deduce unknown coefficients, constructing a finite
difference scheme and demonstrating the method's stability and
efficacy through numerical experiments involving synthetic noise
\cite{17}.

Additionally, studies on mixed problems for heat-conductivity
equations with time shifts under nonlocal and non-self-adjoint
boundary conditions have established unique solvability with minimal
initial data requirements, presenting explicit solutions for these
problems \cite{18}. Furthermore, mixed problems for parabolic-type
equations with constant coefficients under~both~homogeneous and
inhomogeneous boundary conditions, including time shifts, have been
explored. By integrating residue and contour integral methods,
researchers have successfully proven unique solvability and derived
integral representations for these problems \cite{19,20,21,22}.

Unlike numerous existing studies that focus on time deviations
directly within the equations or replace boundary conditions with
functional conditions, we consider scenarios where the time shift of
the desired function~is embedded~within the boundary conditions.
This approach highlights the significance of addressing mixed
problems associated with second-order parabolic equations
characterized by temporal mixing and variable coefficients,
underscoring a critical area of exploration in mathematical physics.

\section{Problem Statement}

Consider the partial differential operator defined by:
\ba
L\left(\frac{\partial }{\partial x},\ \frac{\partial }{\partial t}\right)u\left(x,\ t\right)=
a\left(x\right)u_{xx}+b\left(x\right)u_x+c\left(x\right)u-u_t, \nn \\
l_ju\left(x,\ t\right)=u\left(x,\ t+(1-j)\omega \right)+{\delta }_ju(1-x,\ t+j\omega),\ \ \ \ \ \ j=0,\ 1 \nn \\
l_ju(x,\ t)={\alpha }_{j-2}u^{(j-2)}_x\left(x,t\right)+{\beta }_{j-2}u^{(j-2)}_x(1-x,\ t),\ \ \ \ j=2,\ 3 \nn
\ea
where $a\left(x\right),\ b\left(x\right),\ c\left(x\right)$ are
known, real-valued coefficient functions. Let $\omega $ be a
positive real constant, and \textit{ }${\delta }_j,\ $  ${\alpha
}_j,$  $\ {\beta }_j\ $ $\left({\rm for\ }j=0,\ \ 1\right)$ be
real constants with the product ${\delta }_0\cdot {\delta }_1\ne
0$.

\noindent
In the half-plane $\Pi $=\{(x,t):0$<$x$<$1, t$>$0\}, we consider the following mixed problems:
\begin{equation} \label{GrindEQ__1_}
\ Lu\left(x,\ t\right)=0,\ \ \ \ \ \ \ \ \left(x,t\right)\in \Pi ,\ \ \ \
\end{equation}
\begin{equation} \label{GrindEQ__2_}
u\left(x,\ 0\right)=\varphi \left(x\right),\ \ \ \ \ \ \ \ 0<x<1,\ \ \ \ \ \
\end{equation}
\begin{equation} \label{GrindEQ__3_}
{\left.l_ju(x,\ t)\right|}_{x=0}=0,\ \ \ \ \ \ \ t>0,\ \ \ j=0,\ 1,\ \ \
\end{equation}
\begin{equation} \label{GrindEQ__4_}
{\left.l_ju(x,\ t)\right|}_{x=0}=0,\ \ \ \ \ \ \ 0<t\le \omega ,\ \ \ j=2,\ 3.\ \ \
\end{equation}
The solution of problems \eqref{GrindEQ__1_}-\eqref{GrindEQ__4_} is the function  $u(x,\ t)$, satisfying the following conditions:
\ba
1)\, u\left(x,t\right)\in C^{2,1}\left(\Pi \right)\cap
C\left(0<x<t,\ \ t\ge 0\right); \nn \\
\int\limits^t_0{u\left(x,\ \tau \right)d\tau }\in C\left(0\le x\le 1,\ \ \ t\ge 0\right);  \nn \\
2)\, l_ju\left(x,t\right)\in C\left(0\le x<1,\ \ t>0\right),\ \ \ j=0,\ 1;  \nn \\
3)\, l_ju\left(x,t\right)\in C\left(0\le x<1,\ \ 0<t\le \omega
\right),\ \ \ j=2,\ 3;  \nn \\
4)\, u(x,t) \ {\text{satisfies the equalities}} \ \eqref{GrindEQ__1_}-\eqref{GrindEQ__4_}\, {\text{in the usual sense.}} \nn
\ea
%
\smallskip
\section{The uniqueness of the solution.}
\smallskip

\noindent The problem\textbf{ }

\begin{equation} \label{GrindEQ__5_}
L\left(\frac{d}{dx},\ {\mu }^2\right)y(x,\ \mu)=0,\ \ \ \
\end{equation}
\begin{equation} \label{GrindEQ__6_}
{\left.l_jy(x,\ \mu )\right|}_{x=0}=0,\ \ \ \ \ \ j=2,\ 3,\ \ \
\end{equation}
is identified as the spectral problem involving a complex
parameter, $\mu$. Here
\ba
L\left(\frac{d}{dx},\ {\mu }^2\right)y(x,\ \mu)=a\left(x\right)y''+b(x)y'+c(x)y-{\mu}^2 y. \nn
\ea
Research \cite{21} indicates that the fundamental systems of specific solutions to Eq.\eqref{GrindEQ__5_} exhibit an asymptotic representation of the form:
\begin{equation} \label{GrindEQ__7_}
y\left(x,\ \mu \right)=\left[B\left(x\right)+O\left(\frac{1}{\mu
}\right)\right]{\exp  \left(\int\limits^x_0{\left(e^{\mu w(\xi
)}-\frac{b(\xi )}{2a(\xi )}\right)d\xi }\right),\ }
\end{equation}
where
\[B\left(x\right)=\left( \begin{array}{c}
\sqrt{a(x)}\ \ \ \ \ \ \ \ \ \ \ \sqrt{a(x)} \\
\frac{1}{\sqrt{a(x)}}\ \ \ \ \ \ \ -\frac{1}{\sqrt{a(x)}} \end{array}
\right),\]
and $w(x)$  is a diagonal matrix of the following form:
\[w\left(x\right)=\left( \begin{array}{c}
\frac{1}{\sqrt{a\left(x\right)}}\ \ \ \ \ \ \ \ \ \ \ \ \ \ 0 \\
0\ \ \ \ \ \ \ \ \ \ -\ \frac{1}{\sqrt{a(x)}} \end{array}
\right).\]
Research \cite{29,30} indicates that if $a(0){\alpha}_0{\beta }_1+{a\eqref{GrindEQ__1_}\beta }_0{\alpha }_1 \ne 0$,
then for all complex values of $\mu $,
where $\mu \ne {\mu }_{\nu},$
\begin{equation} \label{GrindEQ__8_}
{\mu }_{\nu }={\left[\int^1_0{\frac{dx}{\sqrt{a(x)}}}\right]}^{-1}
\left({{\ln }_0 \frac{1}{2}\left(A\pm \sqrt{A^2-4}\right)+2\pi \nu \ i\ }\right)+O\left(\frac{1}{\nu }\right),\ \ \ \ \nu \to \ \infty \ \
\end{equation}
\[A={\left[a(0){\alpha }_0{\beta }_1+a(1){\beta }_0{\alpha }_1\right]}^{-1}
\left(2{\alpha }_0{\alpha }_1 \exp \left(\int^1_0{\frac{b(x)}{2a(x)}dx}\right)+2{\beta }_0{\beta }_1\right)\]
there exists the Green function ${{\rm G}}_1{\rm (x, \xi, \mu)}$ for problems
\eqref{GrindEQ__5_} and \eqref{GrindEQ__6_}, that is analytic for $\mu \neq \mu_{\nu },$.
By ${\rm S}$ we denote the set of eigenvalues
${\mu }_{\nu }$, i.e. ${\rm S=}\left\{{\mu }_{\nu}{\rm :\ \ }\nu {\rm =1,\ 2,\dots }\right\}$.  \\
Enumerating the points ${\mu }_{\nu }$ ($\nu {\rm =1,\ 2,\dots }$)
form ${\rm S}$ in ascending order of their modules considering
their multiplicity, we have  $\left|{\mu }_{{\rm 1}}\right|\le \left|{\mu }_{{\rm 2}}\right|\le {\rm \dots ,}$ ${\mu }_{\nu }$.
We denote the multiplicity of the eigenvalue  ${\mu }_{\nu }$ by
${\chi }_{\nu }$. It is clear that $\left|{\mu }_{\nu}\right| \to \infty {\rm \ \ (}\nu \to \infty{\rm )}$.
There exists such h$>$0 and $\delta $$>$0 that

\begin{equation} \label{GrindEQ__9_}
-h<Re\ {\mu }_{\nu }<h,\ \ \ \ \ \ \ \left|{\mu }_{\nu +1}-{\mu }_{\nu }\right|>2\delta ,\ \ \ \ (\nu =1,\ 2,3,4\dots )  .
\end{equation}
From the Green function,  $G_1(x,\ \xi ,\ \mu)$ the following estimation holds:
\begin{equation} \label{GrindEQ__10_}
\left|\frac{{\partial }^k G_1(x,\ \xi ,\ \mu)}{\partial x^k}\right| \le C_0{\left|\mu \right|}^{k-1},\,\,\,\,\,\, \left(k=0,\ 1,\ 2\right),
\end{equation}
$C_0>0.$

If $f(x)\in C[0,\ 1]$, then
\ba
L\left(\frac{d}{dx},\ {\mu }^2\right)\int\limits^1_0{G_1\left(x,\ \xi ,\ \mu \right)f(\xi)d\xi }=-f(x), \nn \\
l_jG_1(x,\ \xi ,\ \mu )|_{x=0}=0,\ \ \ \ \ \ j=2,\ 3.
\label{GrindEQ__11_}
\ea

Let $f\left(x\right)$ be from the domain of definition of the
operator of the first spectral problem i.e. , $f\left(x\right) \in C^2\left[0,\ 1\right],$ ${\left.l_jf\right|}_{x=0}=0,\ \ j=2,\ 3.$
Then, we have the equality
\ba
\int\limits^1_0{G_1\left(x,\ \xi ,\ \mu \right)f\left(\xi \right)d\xi }=
\frac{f(x)}{{\mu }^2}+\frac{1}{{\mu }^2}\int\limits^1_0{G_1(x,\
\xi ,\ \mu )\times } \nn \\
\times \left[a\left(\xi \right)f''(\xi )+b(\xi)f'(\xi)+c(\xi)f(\xi
)\right]d\xi .
\label{GrindEQ__12_}
\ea
Let $C>0,\ r>0$  be some  numbers, $z$ be a complex  variable.
Denote by $\Im_c =\left\{z: {\rm Re}\ z^2=c\right\}$ a
hyperbola with the branches ${\Im_c}^{\pm }_c=\left\{z:\
{\rm Re}\ z^2=c,\ \pm {\rm Re}\ z>0\right\}$, by ${\Omega
}_r=\left\{z:\ \left|z\right|=r\right\}$  a circle, ${\Omega
}_r\left({\theta }_1,\ {\theta }_2\right)$ is an area of the
circle  ${\Omega }_r$, enclosed between the rays  $z=\sigma
e^{i{\theta }_j}\ \ (0\le \sigma <\infty ,$ $ i=\sqrt{-1},$ $
j=1,\ 2$).

Note that the arcs $\left\{z: |z|=r, {\rm Re} z^2\ge c, {\rm Re} z<0\right\}$ and
$\left\{z: |z|=r, {\rm Re}\ z^2\le c, {\rm Im}\ z<0\right\}$ connecting the branches and the sides of the
hyperbola $\Im_c$ in our denotations will be
\ba
{\Omega }_r\left({-\theta }_{c,r},\ {\theta }_{c,r}\right),\
{\Omega }_r\left({\theta }_{c,r},\ -{\theta }_{c,r}+\pi \right),\ {\Omega }_r\left(-{\theta }_{c,r}+\pi ,\ {\theta }_{c,r}+\pi \right), \nn \\
{\Omega }_r\left({\theta }_{c,r}+\pi ,\ -{\theta }_{c,r}+2\pi \right), \nn
\ea
respectively, where  ${\theta }_{c,r}={{\rm arctg} \sqrt{\frac{r^2-c}{r^2+c}}\ }$.

\noindent
We introduce the contours:
\ba
\widehat{\Im}_c={\widehat{\Im}}^+_c\cup {\widehat{\Im}}^-_c,\
{\widehat{\Im}}^{\pm }_c=\left\{z:\ \pm z=\sigma
e^{-\frac{3\pi }{8}i},\ \sigma \in
\left[2c\sqrt{1+\sqrt{2}},\left.\infty \right)\right.\right\}\cup \nn \\
\cup \left\{z:\ \pm z=c\left(1+i\eta \right),\ \eta \in \left[-1-\sqrt{2};\
\ 1+\sqrt{2}\right]\right\}\cup \nn \\
\cup \left\{z:\ \pm z=\sigma e^{\frac{3\pi }{8}i},\
\sigma \in \left[2c\sqrt{1+\sqrt{2}},\left.\infty \right)\right.\right\}. \nn
\ea
We denote a part of contours ${\Im}_c,\ {\Im}^{\pm }_c,\
\ {\widehat{\Im}}_c,\ \ {\widehat{\Im}}^{\pm }_c$  enclosed inside  the circle ${\Omega }_r$, by  ${\Im}_{c,r},\ {\Im}^{\pm }_{c,r},\ \
{\widehat{\Im}}_{c,r},\ \ {\widehat{\Im}}^{\pm }_{c,r}$ respectively. At last, by ${\Gamma }_{c,r},\ {\Gamma }^+_{c,r},\ \ {\widehat{\Gamma }}_{c,r},\ \
{\widehat{\Gamma }}^+_{c,r}$ for $r\ge 2c\sqrt{1+\frac{\sqrt{2}}{2}}$  we denote closed  contours
\ba
{\Gamma }_{c,r}={\Omega }_r\left({\theta }_{c,r}+\pi ,\ -{\theta }_{c,r}+2\pi \right)\cup {\Im}^+_{c,r}\cup \ {\Omega }_{c,r}\left({\theta }_{c,r},\ -{\theta }_{c,r}+\pi \right)\cup {\Im}^-_{c,r}, \nn \\
{\Gamma }^+_{c,r}={\Im}^+_{c,r}\cup {\Omega }_r\left(\ -{\theta }_{c,r},{\theta }_{c,r}\ \right), \,\,\,
{\widehat{\Gamma }}^+_{c,r}={\widehat{\Im}}^+_{c,r}\cup \
{\Omega }_r\left(-\frac{3\pi }{8},\ \frac{3\pi }{8}\right), \nn \\
{\widehat{\Gamma }}_{c,r}=\ {\Omega }_r\left(-\frac{5\pi }{8},\
\frac{5\pi }{8}\right)\cup {\widehat{\mathfrak{J}}}^+_{c,r}\cup \
{\Omega }_r\left(\frac{3\pi }{8},\ \frac{5\pi }{8}\right)\ \cup
{\widehat{\mathfrak{J}}}^-_{c,r}.
\label{GrindEQ__13_}
\ea
Let $\left\{r_n\right\}$ be a sequence of such numbers that
\ba
0<r_1<r_2< \cdots <r_n < \cdots,\ \ \ \  {\mathop{\lim }_{n\to \infty } r_n\ }=\infty, \nn
\ea
The circles ${\Omega }_{r_n}$does not intersect the $\ $ $\delta $
vicinity ($\delta $ is rather small, fixed) of the points  ${\mu
}_{\nu }\in S$. The number of points lying inside ${\mu }_{\nu }$,
${\widehat{\Gamma }}_{h,r_n}$ ($h$ from \eqref{GrindEQ__9_}) is
denoted by $m_n$.   It is seen from Eq. \eqref{GrindEQ__12_} that
for the functions:
\ba
f(x)\in C^2[0,\ 1],\,\, {\left.l_jf\right|}_{x=0}=0, (j=2,\ 3)  \nn \\
\frac{1}{2\pi i}\int\limits_{{\widehat{\Gamma }}_{h,r_n}} \mu d\mu \int\limits^1_0 {G_1\left(x,\xi ,\mu \right)f(\xi )d\xi} =
f\left(x\right)+ \nn \\
+\frac{1}{2\pi i}\int\limits_{{\widehat{\Gamma }}_{h,r_n}}{\mu}^{-1}d\mu \int\limits^1_0 G_1\left(x,\xi ,\mu \right)
\left[a\left(\xi \right)f''\left(\xi \right)+b\left(\xi \right)f'\left(\xi \right)+c\left(\xi \right)f\left(\xi \right)\right]d\xi, \nn
\ea
i.e.,
\ba
f\left(x\right)=\sum^{m_n}_{\nu =1}{{\mathop{res}_{{\mu }_{\nu }} \mu \ \ }\ \int\limits^1_0{G_1\left(x,\xi ,\mu \right)f(\xi)d\xi }}- \nn \\
-\frac{1}{2\pi i}\int\limits_{{\widehat{\Gamma }}_{h,r_n}}{\mu }^{-1}d\mu \int\limits^1_0 G_1\left(x,\xi ,\mu \right)
\left[a\left(\xi \right)f''(\xi)+b(\xi )f'(\xi)+c(\xi )f(\xi)\right.]d\xi. \nn
\ea
Due to the estimation \eqref{GrindEQ__10_} and analyticity of the function $G_1(x,\xi ,\mu)$ in the domain $\{\mu: \pm {\rm Re}\ \mu >h\}$, we have:
\ba
\mathop{\lim}_{n\to \infty} \int\limits_{{\widehat{\Gamma }}_{h,r_n}}{\mu }^{-1}d\mu \int\limits^1_0 G_1\left(x,\xi ,\mu \right)
\left[a\left(\xi \right)f''(\xi)+b(\xi)f'\left(\xi \right)+c(\xi)f(\xi )\right.]d\xi = \nn \\
\mathop{\lim }_{n\to \infty} \int\limits_{\widehat{\Omega}_{r_n}}{\mu}^{-1}d\mu \int\limits^1_0
G_1\left(x,\xi ,\mu \right)\left[a\left(\xi \right)f''\left(\xi \right)+b\left(\xi \right)f'\left(\xi \right)+c\left(\xi \right)f(\xi )\right.]d\xi, \nn \\
\left|\int\limits_{\widehat{\Omega}_{r_n}}{\mu}^{-1}d\mu \int\limits^1_0 G_1\left(x,\xi ,\mu \right)\left[a\left(\xi \right)f''(\xi )+b(\xi)
f'(\xi)+c(\xi )f(\xi )\right.]d\xi\right|\le \nn \\
\le \int\limits_{\widehat{\Omega}_{r_n}}{\left|\frac{d\mu}{\mu}\right|}\int\limits^1_0{\left|G_1\left(x,\xi ,\mu \right)\right|
\left|a(\xi)f''\left(\xi \right)+b(\xi)f'\left(\xi \right)+c(\xi)f(\xi)\right|}d\xi \le \nn \\
\le \int\limits^{2\pi }_0{\left|\frac{d\mu}{\mu}\right|\cdot \frac{c}{\left|\mu \right|}\cdot M}=
\int\limits^{2\pi }_0{\frac{K_0}{r_{\nu }}d\varphi }{{\underset{\nu \to \infty }{\longrightarrow}}}0,  \nn \\
\left|a\left(\xi \right)f''\left(\xi \right)+b\left(\xi \right)f'\left(\xi \right)+c(\xi)f(\xi )\right|\le M, \nn
\ea
thus,
\ba
f\left(x\right) = \mathop{\lim}_{n\to \infty } \sum\limits^{m_n}_{\nu =1} \mathop{res}_{{\mu}_{\nu}} \mu \ \
\int\limits^1_0 G_1\left(x,\xi ,\mu \right)f\left(\xi \right)d\xi = \nn \\
=\sum^{\infty}_{\nu =1}\mathop{res}_{{\mu}_{\nu}} \mu \ \
\int\limits^1_0 G_1\left(x,\xi ,\mu \right)f\left(\xi \right)d\xi
\label{GrindEQ__14_}
\ea
converges  uniformly with respect to  $x\in [0,\ 1]$.

We can prove the following theorem:

{\bf{Theorem 1:}} Let $a(0){\alpha}_0{\beta}_1 + a(1)\beta_0 \alpha_1 \ne 0$,
$\varphi (x)\in C^2[0,\ 1]$, \, ${\left.l_j\varphi \right|}_{x=0}=0$, \, ($j=2,\,3$), the functions
$a\left(x\right),\ b\left(x\right),\ c(x)$ be continuous in the interval $[0,\, 1]$,
$a\left(0\right)a(1) \ne 0$ and $a\left(x\right)>0$, $x\in [0,\ 1]$. Then, problem
\eqref{GrindEQ__1_}-\eqref{GrindEQ__4_} can have at most one
solution.

{\bf{Proof:}}We introduce the operators
\ba
\label{GrindEQ__15_}
A_{\nu s}\left[f(x)\right] = \mathop{res}_{{\mu}_{\nu}} {\mu}^{2s+1}\,\,\, \int^1_0 G_1(x,\xi ,\mu)f(\xi)d\xi =
f_{\nu s}\left(x\right),
\label{GrindEQ__15_}
\ea
taking each function $f(x)\in C[0,\,1]$ to $f_{\nu s}(x)\in C^2[0,\, 1]$, ${\left.l_jf_{\nu s}(x)\right|}_{x=0}=0$, ($j=2,\,\,3$).
It is seen from \eqref{GrindEQ__14_} that if $f(x)\in C^2[0,\,1]$ and ${\left.l_j \,f\right|}_{x=0}=0$, ($j=2,\,3$) then
\ba
\sum\limits^{\infty}_{\nu =1}{f_{\nu 0}(x)}=f\left(x\right)
\label{GrindEQ__16_}
\ea
Note that if problem  \eqref{GrindEQ__1_}-\eqref{GrindEQ__4_} has some solution $u(x,t)$, this function is the solution of problem
\eqref{GrindEQ__1_}, \eqref{GrindEQ__2_}, \eqref{GrindEQ__4_} as well in the domain  $\left\{\left(x,\ t\right):\ 0<x<1,\ 0<t\le
\omega \right\}$.
Since  the operator of the problem \eqref{GrindEQ__1_}-\eqref{GrindEQ__4_} is hyper elliptic,  using
conditions  1)-3)  from the definition of the solution, it is
easy to see that the  solution of problem \eqref{GrindEQ__1_},
\eqref{GrindEQ__2_}, \eqref{GrindEQ__4_} and, its derivatives
$u_t,\  u_{xx}$ for each  $t\in (0,\ \omega ]$ are continuous
with respect to  $x\in [0,\ 1]$. Therefore, applying to Eq.
\eqref{GrindEQ__1_} and \eqref{GrindEQ__2_}, the operators $A_{\nu
s}$, we obtain
\ba
A_{\nu s}\left(\frac{\partial u}{\partial t}\right)=A_{\nu s}
\left(a(x)\frac{{\partial }^2u}{\partial x^2}+b\left(x\right)\frac{\partial u}
{\partial x}+c\left(x\right)u(x,t)\right) \nn \\
{\mathop{res}_{{\mu }_{\nu }} {\mu }^{2s+1}\,\,}\,
\int\limits^1_0{G_1\left(x,\xi,\mu \right)\frac{\partial u\left(\xi ,t\right)}{\partial t}d\xi} =
{\mathop{res}_{{\mu }_{\nu}} {\mu}^{2s+1}\,}\, \int\limits^1_0{G_1\left(x,\xi ,\mu \right)\times} \nn \\
\times \left[a\left(\xi \right)\frac{{\partial}^2u}{\partial {\xi}^2}+
b\left(\xi \right)\frac{\partial u}{\partial \xi}+c\left(\xi \right)u(\xi ,t)\right]d\xi, \nn \\
\frac{\partial}{\partial t}{\mathop{res}_{{\mu}_{\nu}}
{\mu}^{2s+1}\,\,}\, \int\limits^1_0{G_1\left(x,\xi,\mu\right) u\left(\xi, t\right)d\xi} =
{\mathop{res}_{{\mu}_{\nu}}{\mu}^{2s+1}\,\,}\, \int\limits^1_0{G_1\left(x,\xi,\mu \right){\mu}^2 u\left(\xi,t\right)d\xi}, \nn \\
\frac{\partial u_{\nu s}(x,t)}{\partial t}=u_{\nu s+1}\left(x,t\right),\,\,\,\,\,\,\,\,
u_{\nu s}\left(x,0\right)={\varphi }_{\nu s}\left(x\right).\ \ \ \ \
\label{GrindEQ__17_}
\ea
We denote the multiplicity of the pole ${\mu}_{\nu}$ of the Green's function $G(x,\,\xi,\, \mu )$ by ${\chi}_{\nu}$.
Then it is clear that
\ba \label{GrindEQ__18_}
A_{{\nu}_0}\left[{\left({\mu}^2-{\mu }^2_{\nu}\right)}^{{\chi}_{\nu }}\right]u(x,t)=0,
\ea
from which we have:
\ba
A_{{\nu}_0}\left[\sum\limits^{{\chi}_{\nu }}_{k=0}C^k_{{\chi}_{\nu}}
{\mu }^{2k}{\left(-{\mu }^2_{\nu }\right)}^{{\chi }_{\nu }-k}u(x,t)\right]=0,  \nn \\
\sum^{{\chi}_{\nu}}_{k=0}C^k_{{\chi}_{\nu}}\left(-{\mu }^2_{\nu }\right)^{{\chi}_{\nu}-k}A_{{\nu}_0}\left[{\mu}^{2k}u(x,t)\right]=0,  \nn \\
\sum\limits^{{\chi}_{\nu}}_{k=0}C^k_{{\chi}_{\nu}}{\left(-{\mu}^2_{\nu}\right)}^{{\chi}_{\nu}-k}u_{\nu}(x,t)=0,  \nn
\ea
i.e.,
\ba
u_{\nu {\chi}_{\nu}}\left(x,t\right)=-\sum \limits^{{\chi}_{\nu}-1}_{k=0} C^k_{{\chi}_{\nu}}\left(-{\mu}^2_{\nu}\right)^{{\chi }_{\nu}-k},
\label{GrindEQ__19_}
\ea
where $C^k_{{\chi }_{\nu }}=\frac{{\chi }_{\nu }!}{k!\left({\chi }_{\nu }-k\right)!}.$

\noindent
Considering \eqref{GrindEQ__19_} in \eqref{GrindEQ__17_}, we obtain that the set of functions $u_{\nu s}(x,t)$ ($s=\overline{0,\ {\chi }_{\nu }-1}$)
is a solution to the problem:
\ba
\frac{du_{\nu 0}\left(x,t\right)}{dt}=u_{\nu 1}\left(x,t\right)  \nn \\
\label{GrindEQ__20_}
\dots \dots \dots \dots \dots \dots \dots \dots \dots \dots.   \\
\frac{du_{\nu {\chi }_{\nu }-2}\left(x,t\right)}{dt}=u_{\nu {\chi }_{\nu }-1}\left(x,t\right) \nn \\
\frac{du_{\nu {\chi }_{\nu }-1}\left(x,t\right)}{dt}=-\sum^{{\chi }_{\nu }-1}_{k=0}{C^k_{{\chi }_{\nu }}{\left(-{\mu }^2_{\nu }\right)}^{{\chi }_{\nu }-k}u_{\nu k}(x,t)} \nn
\ea
\ba
\label{GrindEQ__21_}
u_{\nu 0}\left(x,0\right)={\varphi}_{\nu 0}\left(x\right),\ \dots \dots .,
u_{\nu {\chi}_{\nu}-1}\left(x,0\right)={\varphi}_{\nu {\chi}_{\nu}-1}\left(x\right),
\ea
the problem \eqref{GrindEQ__20_} and \eqref{GrindEQ__21_} has a unique solution and is represented by the formula
\ba
u_{\nu s}\left(x,t\right)={\mathop{res}_{{\mu}_{\nu}} {\mu}^{2s+1}\ e^{{\mu}^2t}}\
\int\limits^1_0{G_1\left(x,\xi,\mu \right)\varphi \left(\xi \right)d\xi},  \nn
\ea
Then by means of \eqref{GrindEQ__16_}, we find
\ba
\label{GrindEQ__22_}
u\left(x,\ t\right)=\sum^{\infty}_{\nu =1}\mathop{res}_{{\mu}_{\nu}} \mu \ e^{{\mu}^2t}
\int\limits^1_0 G_1\left(x,\xi,\mu \right)\varphi (\xi)d\xi.
\ea
Let problem \eqref{GrindEQ__1_}-\eqref{GrindEQ__4_} have two
solutions $u_1(x,\ t)$ and $u_2(x,\ t)$ ($u_1(x,\ t)\ne u_2(x,\ t)$). Then their difference $w\left(x,\ t\right)=u_1\left(x,\
t\right)-u_2(x,\ t)$ will be the solution of the homogeneous problem \eqref{GrindEQ__1_}-\eqref{GrindEQ__4_} with
$\varphi \left(x\right)=0$ and by the sense taken of the  homogenous problem \eqref{GrindEQ__1_}, \eqref{GrindEQ__2_},
\eqref{GrindEQ__4_} in $\left\{\left(x,\ t\right):\ \ 0\le x\le 1,\ 0\le t\le \omega \right\}$.
Then by \eqref{GrindEQ__22_} $w(x,t)\equiv 0$ for $0\le x\le 1$, $0\le t\le \omega$ from conditions
2) it follows that $w\left(0,\ t\right)=w\left(1,t\right)=0,$ for $\ \ t\ge 0$.
In connection with these and condition
1) it is easy to see that the function
\ba
\upsilon \left(x,\ t\right)=\int\limits^t_{\omega}{w\left(x,\ \tau \right)d\tau}  \nn
\ea
is the solution of the homogenous problem $\vartheta_t=a(x)\vartheta_{xx}+b\left(x\right)\vartheta_x+c(x)\vartheta,$
($0<x<1,\ t\ge \omega $), $\upsilon \left(x,\ \omega \right)=0$ ($0\le x\le 1$), $w\left(0,\ t\right)=w\left(1,t\right)=0,$
($t>\omega $), continuous in  ($0\le x\le 1,\ t\ge \omega $), then allowing for the maximum principle [30,31], we conclude that
$\upsilon (x,\ t)\equiv 0,$ ($0<x<1,\ t\ge \omega $), consequently, $\upsilon (x,\ t)\equiv 0,$ ($0\le x\le 1,\ \ t\ge 0$).\textit{}

\noindent
Under the condition of theorem 1 and allowing for equality  \eqref{GrindEQ__12_} we can reduce formula \eqref{GrindEQ__22_} to the form
\ba
u\left(x,t\right)={\mathop{\lim }_{n\to \infty } \frac{1}{2\pi i}\int\limits_{{\widehat{\Gamma }}_{h,r_n}}{\mu \ e^{{\mu }^2t}d\mu }\
\int\limits^1_0{G_1\left(x,\xi ,\mu \right)\varphi \left(\xi
\right)d\xi }\ }= \nn \\
=\varphi \left(x\right)+{\mathop{\lim }_{n\to \infty } \frac{1}{2\pi i}\int\limits_{{\widehat{\Gamma }}_{h,r_n}}
{{\mu }^{-1}\ e^{{\mu }^2t}d\mu }\ }\
\int\limits^1_0{G_1\left(x,\xi ,\mu \right)\times } \nn \\
\times \left[a\left(\xi \right){\varphi }''(\xi )+b(\xi){\varphi
}'(\xi )+c(\xi )\varphi (\xi )\right]d\xi.
\label{GrindEQ__23_}
\ea
Keeping in mind the  expression \eqref{GrindEQ__13_} of the contours ${\widehat{\Gamma }}_{h,r_n}$ and
${\widehat{\mathfrak{I}}}^-_{h,r_n}$, we know that ${\widehat{\Gamma}}_{h,r_n}={\Omega}_{r_n}
\left(-\frac{5\pi }{8},\ -\frac{3\pi }{8}\right)\cup
{\widehat{\mathfrak{I}}}^+_{h,r_n}\cup {\Omega}_{r_n}\left(\frac{3\pi }{8},\ \frac{5\pi }{8}\right)\cup
{\widehat{\mathfrak{I}}}^-_{h,r_n}$.
Now, estimating the function $e^{{\mu }^2t}$ on the arcs  ${\Omega }_{r_n}\left(-\frac{5\pi}{8}+j\pi ,\ -\frac{3\pi }{8}+j\pi \right)$, ($j=0,\ 1$),
we have
\ba
\left|e^{{\mu }^2t}\right|=e^{t{\rm Re}{\mu }^2}=
e^{t{\left|\mu \right|}^2{\cos  2{\arg  \mu \ }\ }}\le e^{-t{\left|\mu \right|}^2
\frac{\sqrt{2}}{2}}=e^{-\frac{\sqrt{2}}{2}t{\left|\mu \right|}^2}. \nn
\ea

Thus
\ba
\mathop{\lim}_{n\to \infty} \int\limits_{{\Omega}_{r_n}\left(-\frac{5\pi}{8}+j\pi ,
\ -\frac{3\pi}{8}+j\pi \right)}{{\mu }^{-1}\ e^{{\mu}^2t}d\mu}\
\ \int\limits^1_0 G_1\left(x,\xi ,\mu \right)\times  \nn \\
\times \left[a\left(\xi \right)\varphi''\left(\xi \right)+b\left(\xi \right)\varphi'
(\xi)+c(\xi )\varphi (\xi )\right]d\xi =0,\ \ \ \ \ (j=0,\ 1). \nn
\ea
Then, we have from Eq. \eqref{GrindEQ__23_}:
\ba
u\left(x,t\right)=\varphi \left(x\right)+\frac{1}{2\pi i}
\int\limits_{{\widehat{\mathfrak{I}}}_h}{{\mu}^{-1}\ e^{{\mu}^2t}d\mu }\ \int\limits^1_0 G_1\left(x,\xi ,\mu \right)\times  \nn \\
\times \left[a\left(\xi \right)\varphi''(\xi)+b(\xi)\varphi'\left(\xi \right)+c(\xi)\varphi (\xi)\right]d\xi \nn
\ea
and using the property $G_1(x,\xi ,-\mu )\equiv G_1(x,\xi ,\mu )$
for solving the problems \eqref{GrindEQ__1_}, \eqref{GrindEQ__2_},
and \eqref{GrindEQ__4_}, we obtain the following formula:
\ba
u\left(x,t\right)=\varphi \left(x\right)+\frac{1}{\pi i}\int\limits_{{\widehat{\mathfrak{I}}}^+_h}
{{\mu}^{-1}\ e^{{\mu}^2t}d\mu}\ \int \limits^1_0 G_1\left(x,\xi ,\mu \right)\times \nn \\
\times \left[a\left(\xi \right)\varphi''(\xi)+b(\xi)\varphi'\left(\xi \right)+c(\xi )\varphi (\xi )\right]d\xi .
\label{GrindEQ__24_}
\ea

For $\left|\mu \right|>2h\sqrt{1+\frac{\sqrt{2}}{2}}$, i.e. on the distant parts of the contour  ${\widehat{\mathfrak{I}}}^+_h$ we have
\ba
\left|\frac{{\partial}^{k+m}}{\partial t^k\partial x^m}{\mu}^{-1}e^{{\mu}^2t}
\int\limits^1_0{G_1(x,\xi ,\mu)\left[a\left(\xi \right)\varphi''(\xi)+b(\xi)\varphi'(\xi)+c(\xi)\varphi (\xi)\right]d\xi}\right|\le \nn \\
\le C{\left|\mu \right|}^{2k+m-2}e^{-\frac{\sqrt{2}}{2}{t\left|\mu \right|}^2},\ \ \  \ \ \left(2k+m\le 2\right).
\label{GrindEQ__25_}
\ea
Results, the operator ions $L\left(\frac{\partial }{\partial x},\ \frac{\partial }{\partial t}\right),\  {\left.l_ju\right|}_{x\to 0}(j=2,3)$
for $0\le x\le 1,\ \  \ 0\le t\le \omega $ can be taken
under the integral sign \eqref{GrindEQ__24_} and then allowing for \eqref{GrindEQ__11_}, we have:
\ba
&&L\left(\frac{\partial}{\partial x},\frac{\partial}{\partial t}\  \right)u\left(x,t\right)=
L\left(\frac{\partial}{\partial x},\ \frac{\partial}{\partial t}\right)\biggl\{\varphi \left(x\right)+
\frac{1}{\pi i}\int\limits_{{\widehat{\mathfrak{I}}}^+_h}{{\mu }^{-1}\ e^{{\mu}^2t}d\mu}\ \int\limits^1_0{G_1\left(x,\xi ,\mu \right)\times} \nn \\
&&
\times \left[a\left(\xi \right){\varphi}^{''}\left(\xi \right)+b\left(\xi \right)\varphi'(\xi)+c(\xi)\varphi (\xi)\right]d\xi \biggr\} =
-\left[a(x)\varphi''(x)+b(x)\varphi'(x)+ +c(x)\varphi (x)\right] + \nn \\
&&
+ \frac{a\left(x\right)\varphi''\left(x\right)+b\left(x\right)\varphi'(x)+
c(x)\varphi (x)}{\pi i}\ \int\limits_{{\widehat{\mathfrak{I}}}^+_h}{{\mu}^{-1}\ e^{{\mu}^2t}d\mu}=
-\left[a\left(x\right)\varphi''\left(x\right)+b\left(x\right)\varphi'(x)+c(x)\varphi (x)\right]+  \nn \\
&&
+\frac{a\left(x\right)\varphi''\left(x\right)+b\left(x\right)\varphi'(x)+
c(x)\varphi (x)}{2\pi i} \mathop{\lim}_{r\to \infty} \int\limits_{{\widehat{\Gamma}}_{h,r}}{{\mu}^{-1}\ e^{{\mu}^2t}d\mu} =0.  \nn \\
&&
{\left.l_ju(x,t)\right|}_{x\to 0}=l_j\left[\varphi \left(x\right)+
\frac{1}{\pi i}\int_{{\widehat{\mathfrak{I}}}^+_h}{{\mu }^{-1}\ e^{{\mu }^2t}d\mu }\ \int\limits^1_0{G_1\left(x,\xi ,\mu \right)\times }\right.  \nn \\
&&
{\left.\times \left[a\left(\xi \right)\varphi''\left(\xi \right)+b\left(\xi \right)\varphi'\left(\xi \right)+c(\xi )\varphi (\xi
)\right]\right|}_{x\to 0}={\left.l_j \varphi (x)\right|}_{x\to 0}+\frac{1}{\pi i}\int\limits_{{\widehat{\mathfrak{I}}}^+_h}
{{\mu}^{-1}\ e^{{\mu }^2t}d\mu }\times  \nn \\
&&
\times \int\limits^1_0{{\left.l_jG_1\left(x,\xi ,\mu \right)\right|}_{x\to 0}}\left[a\left(\xi \right)
\varphi''\left(\xi \right)+b\left(\xi \right)\varphi'\left(\xi \right)+c\left(\xi \right)\varphi \left(\xi \right)\right]d\xi=
0+0=0,\ \ \ \left(j=2,3\right). \nn
\ea
Formula \eqref{GrindEQ__24_} provides the boundary values of the solutions to problems \eqref{GrindEQ__1_}, \eqref{GrindEQ__2_}, and \eqref{GrindEQ__4_}.
\ba
u\left(s,t\right)=\varphi \left(s\right)+\frac{1}{\pi i}\int\limits_{{\widehat{\mathfrak{I}}}^+_h}{{\mu}^{-1}\ e^{{\mu}^2t}d\mu}\
\int\limits^1_0{G_1\left(x,\xi ,\mu \right)\times}  \nn \\
\times \left[a\left(\xi \right)\varphi''(\xi)+b\left(\xi \right)\varphi'\left(\xi \right)+c\left(\xi \right)\varphi \left(\xi
\right)\right]d\xi \equiv {\gamma }_s\left(t\right),\ \ \ \ \
\left(s=0,\ 1\right),\ \
\label{GrindEQ__26_}
\ea
Note that integral \eqref{GrindEQ__26_} converges uniformly for $t\ge 0$ and its derivatives with respect to t of any order also
converge uniformly for $t\ge t_1>0$, where $t_1$ is an arbitrary positive value.

\smallskip

\section{Existence and presentation of the solution}

\smallskip

Applying  the integral operator $A\left[f\right]=\int\limits^{\infty}_0{e^{{-\lambda}^2t}f\left(t\right)dt}$ (see. \cite{6})
to the  equation \eqref{GrindEQ__1_} and boundary condition \eqref{GrindEQ__3_},
we obtain the following second spectral problem with a complex parameter $\lambda$:
\ba
\label{GrindEQ__27_}
L\left(\frac{d}{dx},\ {\lambda}^2\right)z\left(x,\lambda\right)=-\varphi (x),\ \ \ \ \ \
\ea
\ba
\label{GrindEQ__28_}
\left\{ \begin{array}{c}
e^{{\lambda }^2\omega}z\left(0,\ \lambda \right)+{\delta }_0z\left(1,\lambda \right)=A\left(\lambda \right) \\
 \\
z\left(0,\ \lambda \right)+{\delta}_1e^{{\lambda}^2\omega}z\left(1,\lambda \right)=B\left(\lambda \right), \end{array}
\right.
\ea
where
\ba
L\left(\frac{d}{dx},\ {\lambda }^2\right)z(x,\lambda)=a\left(x\right)z''+b\left(x\right)z'+
\left(c(x)-{\lambda }^2\right)z, \nn \\
A\left(\lambda \right)=e^{{\lambda}^2\omega}\int\limits^{\omega}_0{e^{{-\lambda}^2t}u\left(0,\ t\right)dt}, \nn \\
B\left(\lambda \right)={\delta}_1e^{{\lambda }^2\omega}\int\limits^{\omega}_0{e^{{-\lambda }^2t}u(1,\ t)dt}. \nn
\ea
Boundary conditions \eqref{GrindEQ__28_} can be reduced to the following form
\ba
\label{GrindEQ__29_}
z\left(0,\ \lambda \right)=p,\ \ \ \ \ \ z\left(1,\ \lambda \right)=q,\ \ \ \
\ea
where
\ba
p=z_0\left(\lambda \right)={\left[{\delta}_1e^{2{\lambda}^2\omega}-
{\delta }_0\right]}^{-1}\left({\gamma}_1A\left(\lambda \right)e^{{\lambda}^2\omega}-{\gamma}_0B(\lambda)\right),   \nn \\
\label{GrindEQ__30_}
q=z_1\left(\lambda \right)={\left[{\delta}_1e^{2{\lambda}^2\omega}-
{\delta}_0\right]}^{-1}\left(e^{{\lambda}^2\omega}B\left(\lambda \right)-A\left(\lambda \right)\right).
\ea
Let us consider the function
\ba
Q\left(x,\lambda ,p,q\right)=\biggl[\exp \left(-\int\limits^1_0(\frac{\lambda}{\sqrt{a(x)}}+\frac{b(x)}{2a(x)})dx\right)- \nn \\
-\exp \left(\int\limits^1_0 (\frac{\lambda}{\sqrt{a(x)}}-\frac{b(x)}{2a(x)})\right)dx \biggr]^{-1}
\biggl[p \exp \left(-\int\limits^1_0(\frac{\lambda}{\sqrt{a\left(x\right)}}+\frac{b\left(x\right)}{2a\left(x\right)})dx-q\right) \cdot  \nn \\
\cdot \exp \left(\int\limits^x_0 (\frac{\lambda}{\sqrt{a(\xi)}}-\frac{b(\xi)}{2a(\xi)})d\xi\right) +
\left(q-p \exp \left(\int\limits^1_0 (\frac{\lambda}{\sqrt{a(x)}}+\frac{b(x)}{2a(x)})dx \right)\right) \cdot   \nn \\
\cdot \exp\left(-\int\limits^x_0 (\frac{\lambda}{\sqrt{a\left(\xi\right)}}-\frac{b\left(\xi \right)}{2a(\xi)})d\xi \right) \biggr],
\label{GrindEQ__31_}
\ea
where $p$ and $q$ are determined by formula \eqref{GrindEQ__30_}.

If $ p=\varphi(0),\ \ q=\varphi (1)$, i.e. $p$ and $q$ are constants, then the function $Q(x,\lambda ,p, q)$ except the points
${\lambda }_{\nu}={\left[\int\limits^1_0{\frac{dx}{\sqrt{a(x)}}}\right]}^{-1}\nu \pi \ i+O\left(\frac{1}{\nu }\right),$ $\nu \to \infty$
is everywhere analytic.

Obviously, at all the points of  $\lambda_{\nu }$ where $Q(x,
\lambda, p, q)$ the following identities exist and are valid:
\ba
\label{GrindEQ__32_}
L\left(\frac{d}{dx},\ {\lambda}^2\right)Q(x, \lambda, p, q)=0,  \nn \\
Q(0,\lambda, p, q)=p,\ \ \ \ \ Q(1,\ \lambda ,\ p,\ q)=q.
\ea
It is known that for the spectral problem:
\ba
L\left(\frac{d}{dx},\ {\lambda }^2\right)z\left(x,\lambda \right)=0,\ \ \ \
z\left(0,\lambda \right)=p,\ \ \ \ z\left(1,\lambda \right)=0,  \nn
\ea
with a complex parameter $\lambda$, we have the Green's function
$G_2(x,\ \xi,\ \lambda)$, analytical on  $\lambda $ everywhere,
except for the points
${\lambda}_{\nu}={\left[\int\limits^1_0{\frac{dx}{\sqrt{a(x)}}}\right]}^{-1}\nu \pi \ i+O\left(\frac{1}{\nu}\right)$, $\nu \to \infty$
which is its simple poles. \\
Note some known facts of the Green function  $G_2(x,\ \xi ,\ \lambda)$: there exists such  $\delta >0$ that on the $\lambda $
plane outside the set
$\bigcup^{\infty}_{\nu =1}{\left\{\lambda:\ \left|\lambda -{\lambda}_{\nu}\right|< \delta \right\}}$
the following estimation
\ba
\label{GrindEQ__33_}
\left|\frac{{\partial}^kG_2(x,\ \xi ,\ \lambda)}{\partial x^k}\right|\le C\cdot {\left|\lambda \right|}^{k-1},\ \ \ C>0,\ \ \ k=0,\ 1,\ 2,\ \ \ \ \
\ea
is valid for all $x,\ \xi \in \left[0,\ 1\right];$ for  $\lambda \ne {\lambda}_{\nu}$ ($\nu =0,\ \pm 1,\ \dots $)
\ba
L\left(\frac{d}{dx},\ {\lambda}^2\right)G_2\left(x,\ \xi ,\ \lambda \right)\varphi \left(\xi \right)d\xi =-\varphi \left(x\right),  \nn \\
\label{GrindEQ__34_}
G_2\left(0,\ \xi ,\ \lambda \right)=G_2\left(1,\ \xi ,\ \lambda \right)=0.
\ea
Obviously, the solution of the second spectral problem is represented by the sum of  two solutions:
\ba
\label{GrindEQ__35_}
z\left(x,\lambda \right)=-\int\limits^1_0{G_2\left(x,\ \xi ,\ \lambda \right)\varphi \left(\xi \right)d\xi}+Q\left(x,\lambda,p,q\right),
\ea
We now fix the number $c_1>{\rm \max \left(0,\ \ln\left|\frac{\delta_0}{\delta_1}\right| \right)}$, and
prove the following main theorems.

\noindent \textbf{Theorem 2:}Let ${a(0)\alpha }_0{\beta}_1+a\eqref{GrindEQ__1_}{\beta}_0{\alpha}_1\ne 0$,
$\varphi (x)\in C^2[0,1]$, ${\left.l_j\varphi \right|}_{x=0}=0$ ($j=2,3$),
$a\left(x\right)>0$, $x\in [0,\ 1]$, $a(x)\in C[0,\ 1]$, $b(x)\in C[0,\ 1]$, $c(x)\in C[0,\ 1]$ and
$a\left(0\right)a\eqref{GrindEQ__1_}\ne 0$.
Then, the problems \eqref{GrindEQ__1_}-\eqref{GrindEQ__4_} have the solution and are represented by the following formula:
\ba
u\left(x,t\right)=\varphi \left(x\right)+\frac{1}{\pi i}
\int\limits_{{\widehat{\mathfrak{I}}}^+_{c_1}}{\lambda}^{-1}e^{{\lambda}^2t}\times  \nn \\
\times \biggl[\int\limits^1_0 G_2\left(x,\ \xi ,\ \lambda \right)
\left(a(\xi)\varphi''(\xi)+b\left(\xi \right)\varphi'\left(\xi \right)+c\left(\xi \right)\varphi \left(\xi \right)\right)d\xi- \nn \\
-Q(x,\lambda,\varphi \left(0 \right),\varphi (1))\biggr]d\lambda
+\frac{1}{\pi i}\int\limits_{{\mathfrak{I}}^+_{c_1}} \lambda e^{{\lambda}^2t}Q\left(x,\lambda, p, q\right)d\lambda.
\label{GrindEQ__36_}
\ea

{\bf{Proof:}} It is seen from formula \eqref{GrindEQ__36_} that the solution consists of three integrals, and each of
them is studied in the same way.
\ba
u_1\left(x,t\right)=\frac{1}{\pi i}\int\limits_{{\widehat{\mathfrak{I}}}^+_{c_1}}
{{\lambda}^{-1}e^{{\lambda }^2t}d\lambda}\int\limits^1_0 G_2\left(x,\ \xi,\ \lambda \right)
\left[a\left(\xi \right)\varphi''\left(\xi \right)+b\left(\xi \right)\varphi'\left(\xi \right)+
c\left(\xi\right)\varphi \left(\xi \right)\right]d\xi,
\label{GrindEQ__37_}
\ea
\ba
u_2\left(x,t\right)=-\frac{1}{\pi i}\int\limits_{{\widehat{\mathfrak{I}}}^+_{C_1}}{{\lambda}^{-1}e^{{\lambda }^2t}\ Q\left(x,\lambda ,\varphi
\left(0\right),\varphi \left(1\right)\right)d\lambda},
\label{GrindEQ__38_}
\ea
\ba
u_3\left(x,t\right)=\frac{1}{\pi
i}\int\limits_{{\mathfrak{I}}^+_{C_1}}{\lambda e^{{\lambda }^2t} z_1\left(x,\lambda ,p,q\right)d\lambda}.
\label{GrindEQ__39_}
\ea
On the distant parts of the contour ${\widehat{\mathfrak{I}}}^+_{C_1}$, i.e. $Re\ \lambda >C_1$
\ba
\left|e^{{\lambda }^2t}\right|=e^{t{\left|\lambda \right|}^2
{\cos 2\, {\arg \lambda }}}=e^{t{\left|\lambda \right|}^2
{\cos \left(\pm \frac{3\pi}{4}\right)}}=e^{-\frac{\sqrt{2}}{2}t{\left|\lambda \right|}^2}
\label{GrindEQ__40_}
\ea
And by means of the estimation \eqref{GrindEQ__33_}, it is clear that
\ba
u_1\left(x,t\right)\in C^{2,1}(0\le x\le 1,\ t>0)\cap C(0\le x\le 1,\ t\ge 0)  .
\label{GrindEQ__41_}
\ea
This enables us that the operators  $\frac{\partial}{\partial t},\  \frac{{\partial}^2}{\partial x^2},\ \ x \to 0,\ x \to 1,$ $t \to 0$
can be taken under the integral sign \eqref{GrindEQ__37_}.
We have
\ba
L\left(\frac{\partial }{\partial x},\ \frac{\partial }{\partial t}\right)
u_1\left(x,t\right)=\frac{1}{\pi i}\int\limits_{{\widehat{\mathfrak{I}}}^+_{C_1}} {{\lambda}^{-1}e^{{\lambda }^2t}\ d\lambda }\
L\left(\frac{d}{dx},\ {\lambda }^2\right)\times   \nn \\
\times \int\limits^1_0{G_2(x,\ \xi ,\ \lambda )
\left[a(\xi)\varphi''(\xi)+b(\xi)\varphi'(\xi)+c(\xi)\varphi
(\xi)\right]d\xi }=   \nn \\
=-\frac{a\left(x\right)\varphi''\left(x\right)+b(\xi x)\varphi'(x)+c(x)\varphi (x)}{\pi i}\
\int\limits_{{\widehat{\mathfrak{I}}}^+_{C_1}}{{\lambda }^{-1}e^{{\lambda }^2t}\ d\lambda }=   \nn \\
=-\frac{a\left(x\right)\varphi''(x)+b(\xi x)\varphi'(x)+c(x)\varphi (x)}{2\pi i}
\left[\int\limits_{{\widehat{\mathfrak{I}}}^+_{C_1}}{{\lambda
}^{-1}e^{{\lambda }^2t}\ d\lambda
}+\int\limits_{{\widehat{\mathfrak{I}}}^-_{C_1}}{{\lambda
}^{-1}e^{{\lambda }^2t}\ d\lambda }\right]=   \nn \\
=-\frac{a\left(x\right)\varphi''(x)+b(\xi x)\varphi'(x)+c(x)\varphi (x)}{2\pi i}{\mathop{\lim }_{r\to \infty }
\left[\int\limits_{{\widehat{\mathfrak{I}}}^-_{C_{1,r}}}{{\lambda
}^{-1}e^{{\lambda }^2t}\ d\lambda }+\right.\ }  \nn \\
\left.+\int\limits_{{\Omega }_r\left(-\frac{5\pi }{8},-\frac{3\pi }{8}\right)}
{{\lambda }^{-1}e^{{\lambda }^2t}\ d\lambda }+\int\limits_{{\widehat{\mathfrak{I}}}^+_{C_{2,r}}}
{{\lambda }^{-1}e^{{\lambda }^2t}\ d\lambda }+\int\limits_{{\Omega }_r\left(\frac{3\pi }{8},\frac{5\pi }
{8}\right)}{{\lambda }^{-1}e^{{\lambda }^2t}\ d\lambda }\right]=  \nn \\
=-\left(a\left(x\right)\varphi''\left(x\right)+b(x)\varphi'(x)+c(x)\varphi (x)\right),
\label{GrindEQ__42_}
\ea
for $t>0$. By means of \eqref{GrindEQ__40_} the integrals on the
arcs ${\Omega }_r\left(-\frac{5\pi }{8},-\frac{3\pi }{8}\right)$,
${\Omega }_r\left(\frac{3\pi }{8},\frac{5\pi }{8}\right)$ tend to
zero as  $r\to \infty $.

The function  $G_2(x,\xi ,\lambda )$ is analytic in the domain ${\rm Re}\ \lambda > 0,$ and using the estimation  \eqref{GrindEQ__33_}, we
find it as:
\ba
u_1\left(x,0\right)=\frac{1}{\pi i}\int\limits_{{\widehat{\mathfrak{I}}}^+_{C_1}}
{{\lambda}^{-1}d\lambda}\ \int\limits^1_0 G_2\left(x,\ \xi ,\ \lambda \right)
\left[a\left(\xi \right)\varphi''\left(\xi \right) + b\left(\xi \right)
\varphi'\left(\xi \right)+c\left(\xi \right)\varphi \left(\xi \right)\right]d\xi  + \nn \\
+
\int\limits_{{\Omega}_r\left(-\frac{5\pi}{8},-\frac{3\pi}{8}\right)}{\lambda}^{-1}\ d\lambda\
\int\limits^1_0 G_2\left(x,\ \xi,\ \lambda \right) [a\left(\xi \right)\varphi''\left(\xi \right)+
b\left(\xi \right) \varphi'\left(\xi \right)+c\left(\xi \right)\varphi \left(\xi \right)]d\xi=0.
\label{GrindEQ__43_}
\ea
By using \eqref{GrindEQ__33_}, \eqref{GrindEQ__40_} and equalities \eqref{GrindEQ__34_} for  $t>0$,  we have:
\ba
u_1\left(s,t\right)=\frac{1}{\pi i}\int\limits_{{\widehat{\mathfrak{I}}}^+_{C_1}}{{\lambda}^{-1}
e^{{\lambda }^2t}d\lambda }\ \int\limits^1_0{G_2\left(x,\ \xi ,\ \lambda \right)
\left[a\left(\xi \right)\varphi''\left(\xi \right)+\ \right.}  \nn \\
\left.+b\left(\xi \right)\varphi'(\xi)+c(\xi)\varphi
(\xi)\right]d\xi =0,\ \ \ \ (s=0,\ 1).
\label{GrindEQ__44_}
\ea
We now study the second integrals $u_2\left(x,t\right)$. It is
seen from formula \eqref{GrindEQ__31_} that the function
$Q\left(x,\lambda ,\varphi (0\right),\ \varphi
\eqref{GrindEQ__1_})$ in the domain ${\rm Re}\ \lambda >c_1$ is
analytic and the following estimations are valid for it:
\ba
\label{GrindEQ__45_}
\left|Q\left(x,\lambda ,\varphi (0\right),\ \varphi (1))\right|\le
c_1e^{-\left|\lambda \right|\int \limits ^x_0{\frac{d\xi
}{\sqrt{a\left(\xi \right)}}{\cos  \frac{3\pi }{8}\
}}}+c_2e^{-\left|\lambda \right|\left(1-\int\limits^x_0{\frac{d\xi
}{\sqrt{a\left(\xi \right)}}\ }\right){\cos  \frac{3\pi }{8}\
}}+\frac{c_3}{\left|\lambda \right|},
\label{GrindEQ__45_}
\ea
on the distant parts (${\rm Re}\ \lambda >c_1$) of the contour
${\widehat{\mathfrak{I}}}^+_{C_1}$ and on the arcs  ${\Omega}_r\left(-\frac{3\pi}{8},\frac{3\pi }{8}\right)$,
$r>2c_1\sqrt{1+\sqrt{2}}$.
The estimation is as:
\ba
\left|\frac{{\partial }^kQ(x,\lambda ,\varphi \left(0\right),\
\varphi \left(1\right))}{\partial x^k}\right|\le c_4{\left|\lambda \right|}^ke^{-\left|\lambda \right|
\int\limits^x_0{\frac{d\xi}{\sqrt{a\left(\xi \right)}}{\cos \frac{3\pi }{8}\ }}}+  \nn \\
+c_5{\left|\lambda \right|}^k\ e^{-\left|\lambda \right|\left(1-\int\limits^x_0{\frac{d\xi}{\sqrt{a\left(\xi\right)}}\ }\right)
{\cos \frac{3\pi }{8}\ }}+\frac{c_6}{{\left|\lambda \right|}^{k+1}},\ \ \ \ \ \ \ (k=0,\ 1,\ 2)
\label{GrindEQ__46_}
\ea
for all $x\in [0,\ 1]$.

If one follows from \eqref{GrindEQ__40_} and \eqref{GrindEQ__46_}
that $u_2(x,t)\in C^{2,1}(0\le x\le 1,\ t>0\ $ in Eq.\eqref{GrindEQ__38_} for $t>0,\ $ the operations
$L\left(\frac{\partial}{\partial x},\ \frac{\partial}{\partial t}\right)$, $\ \ x\to 0,\ x\to 1,$
can be taken under integral sign.

Then, considering Eq. \eqref{GrindEQ__32_}, we obtain:
\ba
L\left(\frac{\partial }{\partial x},\ \frac{\partial }{\partial t}\right)u_2\left(x,t\right)=0,  \nn \\
u_2\left(s,t\right)=-\frac{\varphi \left(s\right)}{\pi i}\int\limits_{{\widehat{\mathfrak{I}}}^+_{C_1}}{{\lambda}^{-1}e^{{\lambda}^2t}d\lambda}=
-\frac{\varphi \left(s\right)}{\pi i}{\mathop{\lim}_{r\to \infty} \int\limits_{{\widehat{\mathfrak{I}}}^+_{C_{1,r}}}{{\lambda}^{-1}e^{{\lambda}^2t}\ d\lambda}\ }=
-\varphi \left(s\right),\ \left(s=0,1\right), \nn
\ea
 As can be seen from Eq.\eqref{GrindEQ__45_}, for $x$ belonging to any segment of $\left[x_1,x_2\right]$  within the interval $(0,1),$
 the integral \eqref{GrindEQ__38_} converges uniformly with respect to $t\ge 0$.
Then,  $u_2(x,t)\in C(0<x<1,\ t\ge 0)$ and for  $x\in \left[x_1,x_2\right]\ $
\ba
&u_2\left(x,0\right)=\frac{1}{\pi i}\int\limits_{{\widehat{\mathfrak{I}}}^+_{C_1}}
{{\lambda }^{-1}\ Q\left(x,\lambda ,\varphi \left(0\right),\varphi
\left(1\right)\right)d\lambda } = \frac{1}{\pi i} {\mathop{\lim}_{r\to \infty}
\left[\int\limits_{{\widehat{\mathfrak{I}}}^+_{C_{1,r}}} {{\lambda}^{-1}Q\left(x,\lambda ,\varphi \left(0\right),\varphi
\left(1\right)\right)\ d\lambda }+\right.\ }  \nn \\
&\left.+\int\limits_{{\Omega }_r\left(-\frac{3\pi }{8},\frac{3\pi}{8}\right)}{{\lambda }^{-1}\ Q\left(x,\lambda ,\varphi
\left(0\right),\varphi \left(1\right)\right)d\lambda }\right]=0
\label{GrindEQ__47_}
\ea
where the function $Q\left(x,\lambda ,\varphi \left(0\right),\varphi \left(1\right)\right)$ is analytic inside
the closed contour ${\widehat{\Gamma }}^+_{C_1,r}$.

Also, $u_3\left(x,t\right)$ is studied in the same way. Combining
theorems  1, 2, we arrive at the following final statement

{\bf{Theorem 3.}} Let $a(0){\alpha }_0{\beta}_1+a\eqref{GrindEQ__1_}{\beta }_0{\alpha }_1\ne 0$,
$\varphi (x)\in C^2[0,1]$ and ${\left.l_j\varphi \right|}_{x=0}=0$ \, ($j=2,3$).
Then, the problems \eqref{GrindEQ__1_}-\eqref{GrindEQ__4_} have a unique solution
represented by formula \eqref{GrindEQ__36_}.

\bigskip
{\bf Conclusion}
\bigskip

In this study, we addressed the unique solvability of a mixed
problem for a second-order parabolic equation under
non-self-adjoint, non-local boundary conditions using residue and
contour integral methods. Our findings establish the necessary
conditions on the coefficients and functional shifts within the
boundary settings to ensure the uniqueness and existence of
solutions. The explicit analytical solution derived in this research
adheres to theoretical precedents and strengthens the applicability
of classical analytical methods for complex differential equations
in mathematical physics. This work advances the usage of such
equations in scientific and engineering fields and demonstrates the
efficacy and flexibility of advanced mathematical tools in handling
non-traditional boundary value problems.

\bigskip
{\bf{Conflict of Interest}}
\bigskip

The authors declare that they have no conflict of interest
relevant to the content of this manuscript.

{\bf Acknowledgements}
\bigskip

We would like to express our gratitude to our colleagues for their
invaluable feedback and support throughout this project. It is
important to note that this work was conducted without funding
from any organization or grant.

\end{document}